\documentclass[12pt]{amsart}
\usepackage{amssymb,graphics}
\usepackage{amsmath}
\usepackage{latexsym}

\newtheorem{thm}{Theorem}
\newtheorem{lem}[thm]{Lemma}

\newtheorem{cor}[thm]{Corollary}

\newcommand{\defeq}{\stackrel{\rm def}{=}}

\def\({\left(}
\def\){\right)}

\def\no={\,{\,|\!\!\!\!\!=\,\,}}
\def\de{\delta}

\def\ff{\mathbf{f}}
\def\bg{\mathbf{g}}
\def\hh{\mathbf{h}}

\def\de{\delta}
\def\De{\Delta}
\def\Ga{\Gamma}

\def\Z{{\mathbb{Z}}}

\def\no={\,{\,|\!\!\!\!\!=\,\,}}

\def\({\left(}
\def\){\right)}

\def\beq{\begin{eqnarray}}
\def\eeq{\end{eqnarray}}

\begin{document}
\date{\today}

\title[A comparison theorem]{A comparison theorem for $f$-vectors \\ of simplicial polytopes}

\author[Anders Bj\"orner]{Anders Bj\"orner}
\email{bjorner@math.kth.se}

\address{Department of Mathematics\\
         Royal Institute of Technology\\
         S-100~44 Stockholm, Sweden}

\maketitle

{\SMALL DEDICATED TO ROBERT MACPHERSON ON THE OCCASION OF 
HIS 60TH BIRTHDAY
}

\begin{abstract}
Let  $f_i(P)$ denote the
number of $i$-dimensional faces of a convex polytope $P$.
Furthermore, let $S(n,d)$ and $C(n,d)$ denote, respectively, 
the stacked and the cyclic $d$-dimensional
polytopes on $n$ vertices. Our main result is that for every
simplicial $d$-polytope $P$, if
$$f_r(S(n_1,d))\le f_r(P) \le f_r(C(n_2,d))
$$
for some integers $n_1, n_2$ and $r$, then 
$$f_s(S(n_1,d))\le f_s(P) \le f_s(C(n_2,d))
$$
for all $s$ such that $r<s$.

For $r=0$ these inequalities are the well-known
lower and upper bound theorems for simplicial polytopes.

The result is implied by a certain ``comparison theorem'' for $f$-vectors,
formulated in Section 4.
Among its other consequences is a similar lower bound theorem for
centrally-symmetric simplicial polytopes.

\end{abstract}

\section{Introduction}

The following extremal problem and its ramifications  
have a long tradition in the theory of convex
polytopes: among all $d$-dimensional polytopes $P$ with $n$
vertices determine the maximum (or, minimum) of $f_i(P)$.
The answers were given around 1970 by McMullen \cite{mcm} and
Barnette \cite{bar}, who proved that (as had been conjectured)  the upper bound
is attained in all dimensions by the cyclic polytope $C(n,d)$
and the lower bound is attained in all dimensions by the stacked polytope $S(n,d)$.

What if we specify the number of $r$-dimensional faces of $P$, for some
$r > 0$, and pose the analogous extremal problem? The following can be said
in general.


\begin{thm}\label{thm1}
Let $P$ be a $d$-dimensional simplicial polytope.\\
Suppose that
$$f_r(S(n_1,d))\le f_r(P) \le f_r(C(n_2,d))
$$
for some integers $n_1, n_2$ and $0\le r\le d-2$.
Then, 
$$f_s(S(n_1,d))\le f_s(P) \le f_s(C(n_2,d))
$$
for all $s$ such that $r<s<d$.
\end{thm}

For $r=0$ these inequalities are the
lower and upper bound theorems of Barnette and McMullen
\cite{bar}, \cite{mcm}, \cite[Ch. 8]{Zi}. The $s=d-1$ case of the upper bound part is
also known; it is covered by the 
``generalized upper bound theorem'' of Kalai \cite[Theorem 2]{Ka}.

The proof of Theorem \ref{thm1} relies on a comparison theorem for $f$-vectors of
simplicial homology spheres (Theorem \ref{comp} in Section 4) 
together with Stanley's proof of necessity for the $g$-theorem \cite{St2}.
By the same technique we obtain the following.
Here $CS(2n,d)$ denotes
the centrally-symmetric stacked $d$-dimensional
polytopes on $2n$ vertices. 

\begin{thm}\label{thm2}
Let $P$ be a $d$-dimensional centrally-symmetric simplicial polytope.\\
Suppose that
$$f_r(CS(2n,d))\le f_r(P) 
$$
for some integers $n$ and $0\le r\le d-2$.
Then, 
$$f_s(CS(2n,d))\le f_s(P) 
$$
for all $s$ such that $r<s<d$.
\end{thm}

\bigskip

{\bf Acknowledgment.}
This paper was written in response to a question of C. Smyth, who asked
whether the upper bound part of Theorem 1 might be true (personal communication).
Partial results in this direction have also been achieved by
A. Werner and G. M. Ziegler (personal communication).
I am grateful to G. M. Ziegler and to an anonymous referee for helpful comments on a
preliminary version of the paper, and to S. Linusson who spotted and helped correct
an error in the proof of Lemma 3.

\section{Preliminaries}

For the standard notions concerning convex polytopes and simplicial
complexes we refer to the literature, see e.g.  \cite{Zi}. 
In this section we gather some basic definitions and recall some core
results.

The {\em cyclic polytope} $C(n,d)$ is defined and extensively discussed in \cite{Zi}.
The {\em stacked polytope} $S(n,d)$, $n>d$, is obtained from the $d$-simplex
by performing an arbitrary sequence of $n-d-1$ stellar subdivisions of facets.
Similarly, the {\em centrally-symmetric stacked polytope} $CS(2n,d)$, $2n\ge 2d$, 
is obtained from the $d$-dimensional cross-polytope 
by performing an arbitrary sequence of $n-d$ pairs of
centrally-symmetric stellar subdivisions of facets. 
For $n>d+1>3$ the combinatorial types of the resulting polytopes
depend on choices made during the
construction, but their $f$-vectors are well-defined.

Let $\De$ be a $(d-1)$-dimensional simplicial complex, and let $f_i$ be the
number of $i$-dimensional faces of $\De$. The sequence $\ff=(f_0,\dots,f_{d-1})$
is called the {\em $f$-vector} of $\De$. We put $f_{-1}=1$.
The {\em $h$-vector} $\hh=(h_0,\dots,h_d)$ of $\De$ is defined by the equation
\begin{equation*}
\sum_{i=0}^{d} f_{i-1} x^{d-i} = \sum_{i=0}^{d} h_i (x+1)^{d-i}. 
\end{equation*}
From now on we fix the integer $d\ge 3$, and let
$\de=\lfloor\frac{d}{2}\rfloor$. The $g$-{\it vector} of 
$\De$ is the integer sequence
$\bg=\(g_0,g_1,\ldots,g_\de\)$ defined by $g_0=1$  and
$$
g_i=h_i-h_{i-1}, \qquad i=1,\ldots,\de.
$$
The $f$-vector, $h$-vector and $g$-vector of a simplicial $d$-polytope 
are those of its boundary complex.

In the case when $\De$ is a homology sphere 
(or, more generally, a psedomanifold such that the complex itself as well as
the link of every face
has the Euler characteristic of a sphere of the same dimension)
we have the {\em Dehn-Sommerville
equations} $h_i =h_{d-i}$, which show that the $f$-vector of $\De$ is completely
determined by its $g$-vector. The linear relation can be expressed as
a matrix product (see e.g. \cite{Bj2} or \cite[p. 269]{Zi})
\begin{equation*}\label{matrix}
\ff=\bg\cdot M_d,
\end{equation*}
where the $\(\de+1\)\times d$-matrix
$M_d=\(m_{ij}\)$ is defined by
\begin{equation*}
m_{i,j}=\binom{d+1-i}{d-j}-\binom{i}{d-j},\qquad\text{for } 0\le
i\le\de, \, 0\le j\le d-1.
\end{equation*}
Thus, the set of $f$-vectors of homology $(d-1)$-spheres coincides 
with the $g$-vector weighted linear
span of the row vectors of $M_d$.

For instance, we have that

\[
M_{10}=\left( \begin{array}{rrrrrrrrrr}
11&55 &165 &330 &462 &462 & 330&165 &55 &11 \\
 1& 10& 45 &120 &210 &252 &210 &120 &45& 9 \\
 0&  1 & 9 &36 &84 & 126 &126 & 84 & 35 & 7 \\
 0  &0 & 1 & 8 & 28 &56 & 70 & 55 &25&5 \\
 0 & 0 & 0 & 1 &7 & 21 & 34 &31 &  15 &3\\
 0 &  0 & 0  &0 & 1 & 5 &10&10 &5 &1  \\
\end{array}
\right)
\]

\bigskip

\section{Nonnegativity of the $M_d$ matrix}

We need the following technical property of the matrix $M_d$.

\begin{lem}\label{minors}
All $2\times 2$ minors of the matrix $M_d$ are
nonnegative.
\end{lem}
\begin{proof}
For $0\le a<b\le\de$ and $0\le r<s\le d-1$,
let $$
\Phi^{a,b}_{r,s} \defeq m_{a,r}m_{b,s} - m_{a,s}m_{b,r}.
$$
We want to show that $\Phi^{a,b}_{r,s} \ge 0$.
\smallskip

\newcommand{\oor}{{\overline r}}
\newcommand{\oos}{{\overline s}}
\newcommand{\tia}{{\tilde a}}
\newcommand{\tib}{{\tilde b}}

Let $\oor\defeq d-r$, $\oos\defeq d-s$,
$\tia\defeq d+1-a$ and $\tib\defeq d+1-b$. 
Then, by definition
\begin{equation*}
\Phi^{a,b}_{r,s} = \left[ \binom{\tia}{\oor}-\binom{a}{\oor} \right]
\left[ \binom{\tib}{\oos}-\binom{b}{\oos}\right] -
\left[ \binom{\tia}{\oos}-\binom{a}{\oos}\right]
\left[ \binom{\tib}{\oor}-\binom{b}{\oor}\right]
\end{equation*}
Rearranging terms, and letting $B^{p,q}_{t,u}$ denote the binomial determinant
$$B^{p,q}_{t,u} \defeq \det \left(\begin{array}{cc}
\binom{p}{t} & \binom{p}{u} \\
 & \\
\binom{q}{t} & \binom{q}{u}
\end{array} \right)
$$
we can write
\begin{equation}\label{eq5}
\Phi^{a,b}_{r,s} = B^{a, \tib}_{ \oos, \oor}
+ B^{\tib, \tia}_{\oos, \oor}
- B^{a, b}_{ \oos, \oor}
- B^{b, \tia}_{ \oos, \oor}
\end{equation}

{\bf Step 1.} 
Note that
\begin{equation}\label{one}
\det
\left(\begin{array}{cc}
m_{i, t} & m_{i, u} \\
m_{j, t} & m_{j, u}
\end{array} \right) \ge 0
\quad\Leftrightarrow\quad \frac{m_{i, t}}{m_{i, u}} \ge \frac{m_{j, t}}{m_{j, u}},
\end{equation}
if $i<j$, $t<u$ and $m_{j,u}>0$.

An elementary argument based on this observation
shows that it suffices to prove nonnegativity of $\Phi^{a,b}_{r,s}$ for
the special case when $b=a+1$.

({\em Remark:} We could also reduce to the case $s=r+1$; however, this leads to no simplification
in what follows.)


\medskip

{\bf Step 2.} 

In order to show that 
$\Phi^{a,a+1}_{r,s} \ge 0$  we put to use the lattice-path
interpretation of binomial determinants, due to Gessel and Viennot \cite{GV}.
\smallskip

Let $L^{p,q}_{t,u}$ denote the set of pairs 
$(P, Q)$ of vertex-disjoint NE-lattice paths in $\Z^2$,
such that $P$ leads from $(0, -p)$ to $(t, -t)$ 
and $Q$ from $(0, -q)$ to $(u,  -u)$. By a
NE-{\em lattice path} we mean a path taking steps N=$(0, 1)$ to the {\em north}
and steps E=$(1, 0)$ to the {\em east}.

\newcommand{\oot}{{\overline{r+1}}}

The formula of Gessel and Viennot \cite[Theorem 1]{GV} states that
$$B^{p,q}_{t,u} = \# L^{p,q}_{t,u}
$$
Thus, from equation (\ref{eq5}) we have
$$\Phi^{a,a+1}_{r,s} = \# L^{a, \tia -1}_{ \oos,  \oor}
+ \# L^{\tia -1, \tia}_{\oos, \oor}
- \# L^{a, a+1}_{ \oos, \oor}
- \# L^{a+1, \tia}_{ \oos, \oor}
$$

For ease of notation we from now let $ L^{p, q} \defeq
L^{p, q}_{ \oos, \oor}$.
The proof will be concluded by producing an injective mapping
$$\varphi:   L^{a, a+1} \cup L^{a+1, \tia}
\rightarrow  L^{a, \tia-1}
\cup L^{\tia-1, \tia}
$$
The construction of the mapping $\varphi$ proceeds by cases.
\smallskip

\noindent
{\bf Case 1:} 
$(P, Q)\in L^{a, a+1}$. Then $\varphi (P, Q)\in L^{a, \tia-1}$ is constructed by
keeping the path $P$ and extending the path $Q$ by an intitial vertical segment
(a sequence of North steps) so that it begins at the point  $(0, -(\tia-1))$.
\smallskip

\noindent
{\bf Case 2:} 
$(P, Q)\in L^{a+1, \tia}$. 
\smallskip

\noindent
{\bf Subcase 2a:} 
Both $Q$ and $P$ begin with N steps.
Then $\varphi (P, Q)\in L^{a, \tia-1}$ is constructed by removing the first step
from both paths.
\smallskip

\noindent
{\bf Subcase 2b:} 
$Q$ begins with an E step. Then $\varphi (P, Q)\in L^{\tia-1, \tia}$ is constructed by
keeping the path $Q$ and extending the path $P$ by an intitial vertical segment
so that it originates in $(0, -(\tia -1))$.
\smallskip

\noindent
{\bf Subcase 2c:} 
$Q$ begins with an N step, and $P$ begins with an E step.
Then $\varphi (P, Q)\in L^{\tia-1, \tia}$ is constructed
as follows. We may assume that $a\ge \oos$, since otherwise some
binomial coefficients are zero and the situation simplifies.
Thus, the path $P$ begins with a sequence of E steps, say
$k$ of them, followed by a $N$ step. Denoting the
rest of $P$ by $P'$ we can write: $P=E^k N P'$.
Similarly, $Q$ has the factorization  $Q=N REN^v E Q'$, where the two E:s
designate the $k$-th and $(k+1)$-st occurrences of the letter ``E'' in $Q$.
See Figure 1 for the geometric idea.

The integers $k$ and $v$ are determined by the definition of the paths $P$ and $Q$.
Let $h$ be the number of occurrences of the letter  ``N'' in $R$.  
Let
 ${\overline P}$ and ${\overline Q}$ be the paths
$${\overline P}= N^{\tia-a-h-3} ERN^2 P' \mbox{ \;  and \;  }
{\overline Q}= E^{k} N^{v}EN^{h+1} Q', $$
originating in the points $(0,-\tia +1)$ and $(0, -\tia)$, respectively.
A straightforward inspection of the construction shows that these paths are
disjoint. Namely, the lowest point on ${\overline P}$ and the highest point on
${\overline Q}$ with first coordinate $k$ are, respectively, $(k, -a-h-2)$ and
$(k, -\tia +v)$. Their distance is $\tia -a-h-v-2 >0$.
Let $\varphi (P, Q)=({\overline P}, {\overline Q}) \in L^{\tia-1, \tia}$.

This defines the mapping $\varphi$ in all cases. Each case separately is
clearly injective.
That there is no interference among the four cases, and hence that $\varphi$ is injective
globally, is most easily seen from following properties of the construction:
\begin{itemize}
\item $\varphi (P,Q) \in L^{a, \tia -1}$ in cases 1 and 2a
\item $\varphi (P,Q) \in L^{\tia -1, \tia}$ in cases 2b and 2c
\item $(0, -a-1)\in \varphi (Q)$ in cases 1 and 2b
\item $(0, -a-1)\notin \varphi (Q)$ in cases 2a and 2c
\end{itemize}This completes the proof.
\end{proof}

\vspace{.8cm}
\begin{center}\label{paths}
$
\begin{array}{ccc}
\resizebox{!}{2.9in}{\includegraphics{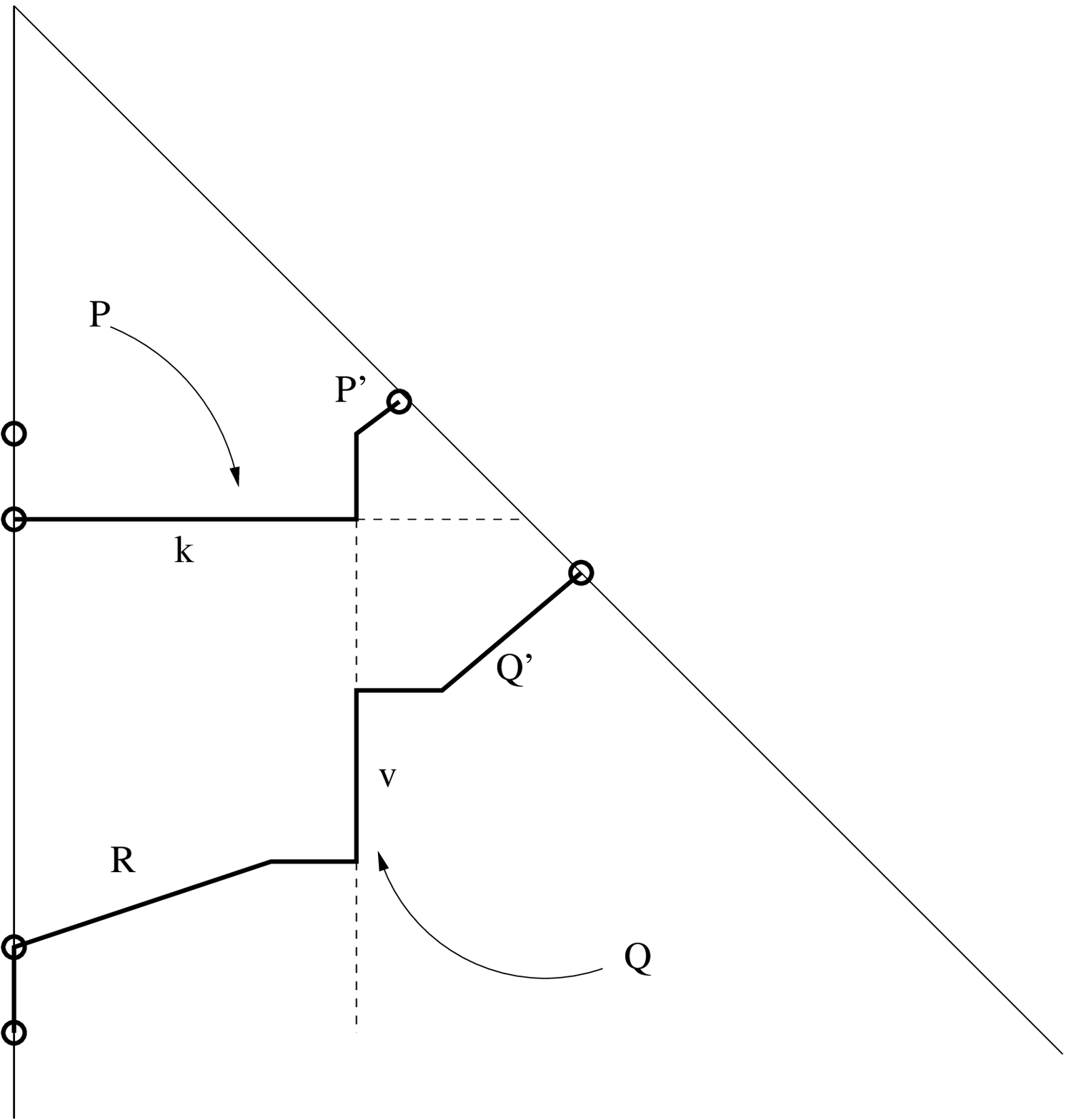}}& & 
\resizebox{!}{2.9in}{\includegraphics{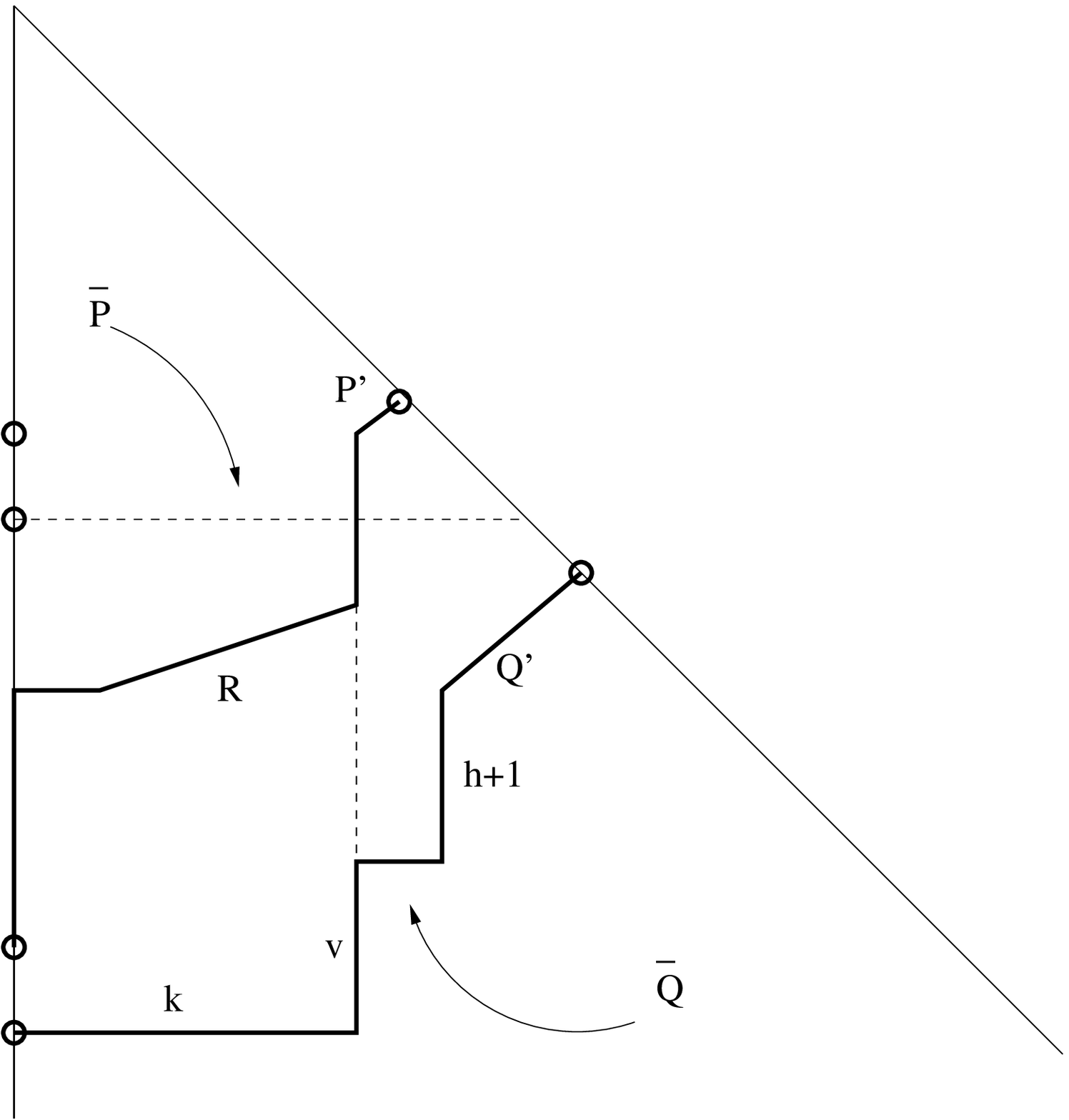}}
\end{array}
$ \\
\vspace{.7cm}
Figure 1: A sketch of subcase 2c.
\end{center}

\medskip

{\bf Remark:} We conjecture that the matrix $M_d$ is {\em totally nonnegative}, meaning
that all minors of all orders are nonnegative. This has been verified 
for all $d\le 13$ by A. Hultman.

\newpage

\section{Homology spheres}

A key role for this paper is played by  the following comparison theorem for $f$-vectors of 
homology spheres. 
 

\begin{thm}\label{comp}
Let $\De$ and $\Ga$ be  $(d-1)$-dimensional simplicial homology spheres 
whose $g$-vectors for some $t$ $(0\le t\le \de)$ satisfy
\begin{itemize}
\item $g_{i} (\De)\ge g_{i}(\Ga)$ for $i=1, \ldots, t$
\item $g_{i} (\De)\le g_{i}(\Ga)$ for $i=t+1, \ldots, \de$.
\end{itemize}
Suppose that
$$f_r(\De)\le f_r(\Ga)
$$
for some $0\le r\le d-2$. Then
$$f_s(\De)\le f_s(\Ga)
$$
for all $s$ such that $r<s<d$.
\end{thm}

\begin{proof}
Let $v_i = g_i (\De)- g_i (\Ga)$. Now,

\beq\label{eq1}
0\ge f_r (\De) -f_r(\Ga) = \sum_{i=0}^{\de} v_i m_{i, r}
= \sum_{i=0}^{\de} v_i m_{i, s} \frac{m_{i, r}}{m_{i, s}}
\eeq
Lemma \ref{minors} implies, in view of equivalence (\ref{one}), that
$$\frac{m_{0, r}}{m_{0, s}} \ge \frac{m_{1, r}}{m_{1, s}} \ge\cdots\ge
\frac{m_{\de, r}}{m_{\de, s}} \ge 0
$$
({\em Remark:} It is possible that $m_{i, s}=0$ for $i=k, \ldots, \de$. Then also
$m_{i, r}=0$ for $i=k-1, \ldots, \de$ while $m_{i, s}>0$ for all $i<v$. This requires
notational adjustments in our argument, but no new ideas.)

By assumption, 
the vector $v=(v_0, v_1, \ldots, v_{\de})$  satisfies 
$$v_1, \ldots, v_{t} \ge 0 \quad\mbox{ and }\quad v_{t+1}, \ldots, v_{\de} \le 0.$$
Thus,
$$  \sum_{i=0}^{\de} v_i m_{i, s} \frac{m_{i, r}}{m_{i, s}} \ge
\left( \sum_{i=0}^{t} v_i m_{i, s}\right) \frac{m_{t, r}}{m_{t, s}} +
\left(\sum_{i=t+1}^{\de} v_i m_{i, s}\right) \frac{m_{t, r}}{m_{t, s}}
$$
which implies that

\begin{equation*}
0\ge f_r (\De) -f_r(\Ga) \ge \frac{m_{t, r}}{m_{t, s}} \left(\sum_{i=0}^{\de} v_i m_{i, s}\right) 
= \frac{m_{t, r}}{m_{t, s}} \left (f_s (\De) -f_s(\Ga) \right)
\end{equation*}
It follows that
$$0\ge f_s (\De) - f_s(\Ga),$$
as desired.
\end{proof}

We will say that an integer vector $(n_0, \ldots, n_{\de})$ is an $m$-sequence
if $n_0 =1$ and $n_j \ge \binom{m}{j}$  implies that
$n_{j-1} \ge \binom{m-1}{j-1}$, for all  $m\ge j>1$. In particular,
if some entry in an $m$-sequence is positive then so are all earlier entries.
The notion of $m$-sequence is less restrictive than the well-established concept
of $M$-sequence, recalled in Section 5.

\begin{cor}\label{cor1} {\rm (Upper bounds)}
Let $\De$ be a $(d-1)$-dimensional homology sphere whose $g$-vector
is an $m$-sequence. Suppose that
$$f_r(\De)\le f_r(C(n,d))
$$
for some integers $n$ and $0\le r\le d-2$. Then
$$f_s(\De)\le f_s(C(n,d))
$$
for all $s$ such that $r<s<d$.
\end{cor}

\begin{proof}
The $g$-vector of the cyclic polytope $C(n,d)$ is
$$g_i (C(n,d)) = \binom{n-d-2+i}{i}
$$
Thus, since $g(\De )$ is an $m$-sequence the conditions of
Theorem \ref{comp} are satisfied.
\end{proof}

Stanley's upper bound theorem for homology spheres \cite{St1} shows that
in the special case when $r=0$ Corollary \ref{cor1} is valid also without the assumption that
$g(\De)$ is an $m$-sequence.

\begin{cor}\label{cor2} {\rm (Lower bounds)}
Let $\Ga$ be a $(d-1)$-dimensional homology sphere whose $g$-vector is
nonnegative.
Suppose that
$$f_r(S(n,d))\le f_r(\Ga)
$$
for some integers $n$ and $r\le d-2$. Then
$$f_s(S(n,d))\le f_s(\Ga)
$$
for all $s$ such that $r<s<d$.
\end{cor}

\begin{proof}
The $g$-vector of the stacked polytope $S(n,d)$ is
$$g_i (S(n,d)) =\left\{
\begin{array}{ll}
1, & \mbox{ for $i=0$} \\
n-d-1,& \mbox{ for $i=1$} \\
0, &\mbox{ for $i>1$}
\end{array}\right. 
$$
Thus, since $g(\Ga)$ is nonnegative the conditions of
Theorem \ref{comp} are satisfied.
\end{proof}

\section{Polytopes}
We recall the definition of an $M$-sequence.
For any integers $k,n\ge 1$ there is a unique way of writing
$$
n=\binom{a_k}{k}+\binom{a_{k-1}}{k-1}+\ldots+\binom{a_i}{i},
$$
so that $a_k>a_{k-1}>\ldots>a_i\ge i\ge 1$. Then define
$$
\partial^k(n)=\binom{a_k-1}{k-1}+\binom{a_{k-1}-1}{k-2}+\ldots
+\binom{a_i-1}{i-1}.
$$
Also let $\partial^k(0)= 0$.

A nonnegative integer sequence $\(n_0,n_1,n_2,\ldots\)$  such that $n_0=1$ and
$$\partial^k\(n_k\)\le n_{k-1}\qquad\text{for
 all } k>1 $$ 
is called an $M$-{\em sequence}.
Clearly, an $M$-sequence is an $m$-sequence (as defined in connection
with Corollary \ref{cor1}), but not conversely.


\medskip

\noindent
{\em Proof of Theorem \ref{thm1}.}\; 
The $g$-vector of a simplicial polytope is an $M$-sequence, by the theorem
of Stanley \cite{St2}. In particular, it is a nonnegative $m$-sequence, so
both Corollaries \ref{cor1} and \ref{cor2} apply.
\hfill$\Box$
\medskip

\noindent
{\em Proof of Theorem \ref{thm2}.}\; 
The $g$-vector of the centrally-symmetric stacked polytope $CS(2n, d)$ is
$$g_i (CS(n,d)) =\left\{
\begin{array}{ll}
1, & \mbox{ for $i=0$} \\
2n-d-1,& \mbox{ for $i=1$} \\
\binom{d}{i} - \binom{d}{i-1}, &\mbox{ for $i>1$}
\end{array}\right. 
$$
Stanley \cite{St3} has shown that
$$g_i (P) \ge \binom{d}{i} - \binom{d}{i-1}, \mbox{ for $i\ge1$}
$$
holds for every centrally-symmetric simplicial polytope $P$.
Hence, Theorem \ref{comp} applies.

\hfill $\Box$


\end{document}